\documentclass[10pt]{amsart}
\usepackage{amsfonts}

\usepackage{amscd}
\usepackage{amsmath}
\usepackage{amsthm}
\usepackage[mathscr]{eucal}

\newtheorem{theo}{Th\'eor\`eme}[section]

\theoremstyle{definition}

\newfont{\got}{eufm10}

\begin{document}
\title{Mesure de Mahler d'hypersurfaces K3}
\author{Marie Jos\' e Bertin}
\date{\today}

\email{bertin@math.jussieu.fr}
\address{Universit\' e Pierre et Marie Curie (Paris 6) }
\maketitle

%%%%%%%%%%%%%%%%%%%%%%%%%%%%%TITRE%%%%%%%%%%%%%%%%%%%%%%%%%% %%%%%%%%%%% 
%MESURE DE MAHLER ET SURFACE K3

%%%%%%%%%%%%%%%%%%%%%%%%%%%%%%%%%%%%%TEXTE%%%%%%%%%%%%%%%%%% %%%%%%%%%%%%%%%%%% 
{\bf R\'esum\'e} Nous exprimons, \` a l'aide de s\'eries d'Eisenstein-Kronecker, la mesure de Mahler de deux familles de polyn\^ {o}mes d\'efinissant des hypersurfaces $K3$ de nombre de Picard g\'en\'erique $19$. Pour certaines de ces surfaces $K3$ singuli\` eres (i.e. de nombre de Picard $20$), nous donnons cette mesure en termes de s\'erie $L$ de Hecke de poids $3$ pour certains {\it Gr\"ossencharacter}.

\section{ Introduction}

 La mesure de Mahler logarithmique $m(P)$ d'un polyn\^{o}me de
Laurent $P\in \mathbb{C}[X_1^{\pm} ,...,X_1^{\pm} ]$est d\'efinie par

$$
m(P)=\frac {1}{(2\pi i)^n} \int_{\mathbb{T} ^n}\log \vert P(x_1^\pm
,...,x_n^\pm)\vert\frac{dx_1}{x_1}...\frac{dx_n}{x_n}
$$
o\`u $\mathbb{T}^n$ d\'esigne le tore $\{(x_1,...,x_n)\in \mathbb{C}^n /
\vert x_1\vert =...=\vert x_n\vert =1\}$. Sa mesure de Mahler $M(P)$ vaut
alors $M(P)=\exp (m(P))$. Si $P$ est un polyn\^{o}me unitaire de $\mathbb{Z}
[X]$, on obtient gr\^{a}ce \`a la formule de Jensen 
$$
M(P)= \prod_{P(\alpha)=0} \max (\vert \alpha \vert , 1),
$$
quantit\'e li\'ee au probl\`eme de Lehmer (1933) sur l'existence d'un
polyn\^{o}me unitaire de $\mathbb{Z} [X]$, non cyclotomique, de mesure de Mahler inf\'erieure
\`a $1,1762...$

 Cependant, gr\^{a}ce \`a un r\'esultat de Boyd \cite{Bo1}\cite{Bo2}, la connaissance
de nombreuses valeurs $M(P)$ pour $P\in \mathbb{Z} [X_1,...,X_n]$ pourrait
\'eclairer le probl\`eme pr\'ec\'edent. 

Depuis quelques ann\'ees, on s'int\'eresse en outre \` a l'obtention
de formules explicites pour $m(P)$ \cite{Bo4}\cite{RV1}\cite{RV2}. Ces formules sont li\'ees \`a
la nature g\'eom\'etrique de la vari\'et\'e alg\'ebrique d\'efinie par $P$.
Par exemple, si $P$ repr\'esente un mod\`ele affine d'une courbe elliptique 
$E$ dont les polyn\^ {o}mes attach\'es aux  faces du polygone de Newton n'ont pour racines que des racines de l'unit\'e et si $P(x,y)\neq 0$ pour tout $(x,y)\in \mathbb{T}^2$, alors $\pi ^2 m(P)
$ est conjectur\'e \^{e}tre un multiple rationnel de la s\'erie $L(E,2)$
associ\'ee \`a la courbe elliptique $E$. Cette conjecture, d\'ecoulant des
conjectures de Beilinson \cite{Bei}, a \'et\'e prouv\'ee par Rodriguez-Villegas dans
certains cas o\`u $E$ poss\`ede de la multiplication complexe \cite{RV1}\cite{RV2}. Elle a
\'et\'e v\'erifi\'ee num\'eriquement par Boyd \cite{Bo3} pour de nombreuses familles
de courbes elliptiques. En outre Rodriguez-Villegas a exprim\'e, pour
certaines familles modulaires de courbes elliptiques, la mesure de Mahler
logarithmique des polyn\^{o}mes associ\'es comme la partie r\'eelle de
certaines s\'eries d'Eisenstein-Kronecker \cite{RV1}\cite{RV2}. 

 Nous nous proposons ici de g\'en\'eraliser ce r\'esultat au cas de
certaines familles de surfaces $K3$ ayant un nombre de Picard g\'en\'erique \'egal \`a $19$
. Nous \'etudierons essentiellement deux familles, la premi\`ere associ\'ee
aux polyn\^{o}mes
$$
P_k=X+\frac {1}{X}+Y+\frac {1}{Y}+Z+\frac {1}{Z}-k
$$
et la seconde li\'ee aux polyn\^{o}mes
 
\begin{multline*}$$
Q_k=X+\frac {1}{X}+Y+\frac {1}{Y}+Z+\frac {1}{Z} \\
+XY+\frac {1}{XY}+ZY+\frac {1}{ZY}+XYZ+\frac {1}{XYZ}-k.$$
\end{multline*}

Nous montrerons les r\'esultats suivants. 

\begin{theo}
 1) Posons $k=t+\frac {1}{t}$ et d\'efinissons 
$$
t=(\frac {\eta (\tau) \eta (6\tau)}{\eta
(2\tau)\eta(3\tau)})^6=q^{1/2}-6q^{3/2}+15q^{5/2}-20q^{7/2}+...
$$
o\`{u} $\eta$ d\'esigne la fonction de Dedekind 
$$
\eta (\tau)=e^{\frac {\pi i \tau}{12}} \prod_{n\geq 1}(1-e^{2\pi i n\tau}).
$$
Alors 
$$
m(P_k)=\Re \{-\pi i \tau+\sum_{n\geq 1}(\sum_{d\mid n}d^3)(\frac
{4q^n}{n}-\frac {16q^{2n}}{2n}+\frac {36q^{3n}}{3n}-\frac {144q^{6n}}{6n})\}.
$$

2) Posons $k=-(t+\frac {1}{t})-2$ et d\'efinissons 
$$
t=\frac {\eta (3\tau)^4\eta (12\tau)^8\eta(2\tau)^{12}}{\eta(\tau)^4\eta
(4\tau)^8\eta (6\tau)^{12}}.
$$
Alors 
$$
m(Q_k)=\Re \{-2\pi i \tau+\sum_{n\geq 1}(\sum_{d\mid n}d^3)(\frac
{-2q^n}{n}+\frac {32q^{2n}}{2n}+\frac {18q^{3n}}{3n}-\frac
{288q^{6n}}{6n})\}.
$$

\end{theo}

\begin{theo}

 Avec les notations du th\'eor\`eme $1.1$, on a l'expression suivante
de la mesure 

 1) 

$$
\begin{aligned}  m(P_k)=&\frac {\Im \tau}{8 \pi ^3}\{ \sum'_{m,\kappa}(-4(2\Re
\frac {1}{(m\tau+\kappa)^3(m\bar{\tau}+\kappa)}+\frac
{1}{(m\tau+\kappa)^2(m\bar{\tau}+\kappa)^2})\\ &+16(2\Re \frac
{1}{(2m\tau+\kappa)^3(2m\bar{\tau}+\kappa)}+\frac
{1}{(2m\tau+\kappa)^2(2m\bar{\tau}+\kappa)^2})\\ &-36(2\Re \frac
{1}{(3m\tau+\kappa)^3(3m\bar{\tau}+\kappa)}+\frac
{1}{(3m\tau+\kappa)^2(3m\bar{\tau}+\kappa)^2})\\ &+144(2\Re \frac
{1}{(6m\tau+\kappa)^3(6m\bar{\tau}+\kappa)}+\frac
{1}{(6m\tau+\kappa)^2(6m\bar{\tau}+\kappa)^2}))\}  \end{aligned} $$

 2) 

$$
\begin{aligned}  m(Q_k)=&\frac {\Im \tau}{8 \pi ^3}\{ \sum'_{m,\kappa}(2(2\Re
\frac {1}{(m\tau+\kappa)^3(m\bar{\tau}+\kappa)}+\frac
{1}{(m\tau+\kappa)^2(m\bar{\tau}+\kappa)^2})\\ &-32(2\Re \frac
{1}{(2m\tau+\kappa)^3(2m\bar{\tau}+\kappa)}+\frac
{1}{(2m\tau+\kappa)^2(2m\bar{\tau}+\kappa)^2})\\ &-18(2\Re \frac
{1}{(3m\tau+\kappa)^3(3m\bar{\tau}+\kappa)}+\frac
{1}{(3m\tau+\kappa)^2(3m\bar{\tau}+\kappa)^2})\\ &+288(2\Re \frac
{1}{(6m\tau+\kappa)^3(6m\bar{\tau}+\kappa)}+\frac
{1}{(6m\tau+\kappa)^2(6m\bar{\tau}+\kappa)^2}))\}  \end{aligned} $$

\end{theo}

Nous terminerons par quelques applications arithm\'etiques.

En particulier, pour certains polyn\^ {o}mes de ces familles d\'efinissant des surfaces $K3$ singuli\` eres (i.e. de nombre de Picard $20$), nous exprimerons leur mesure de Mahler comme des s\'eries $L$ de Hecke pour un {\it Gr\"ossencharacter} de poids $3$.

\section{ Quelques rappels sur les surfaces $K3$}

 Nous allons donner quelques r\'esultats permettant de comprendre les
m\'ethodes utilis\'ees. Le lecteur int\'eress\'e par plus de d\'etails
pourra par exemple consulter  \cite{VY} \cite{Y}. 

 Une surface $K3$ d\'efinie sur $\mathbb{C}$ est une surface 
$X\subset \mathbb{P}^3$ v\'erifiant 
$$H^1(X,\mathcal {O}_X)=0$$
et 
$$K_X=0
$$
(i. e. dont le faisceau canonique est trivial). 

 Une surface $K3$ est dite alg\'ebrique si elle admet un fibr\'e en
droites ample, d\'efinissant un plongement projectif de $X$ dans un espace
projectif. Le caract\`ere alg\'ebrique est caract\'eris\'e par le fait que
le degr\'e de transcendance de son corps de fonctions $\mathbb{C}(X)$ vaut $2
$. Par exemple, un rev\^etement double ramifi\'e le long d'une sextique
plane est une surface $K3$. C'est ainsi le cas de la surface dont un
mod\`ele affine est donn\'e par le polyn\^{o}me $P_k$, pour $k\neq \pm 2,\pm
6$, car il s'\'ecrit 

$$
(2Z+X+\frac {1}{X}+Y+ \frac {1}{Y}-k)^2=(X+\frac {1}{X}+Y+\frac
{1}{Y}-k)^2-4.$$

Si $X$ est une surface $K3$, il existe une unique 2-forme holomorphe $\omega$
sur $X$, unique \`a un facteur scalaire pr\`es.

Par exemple si $F(X_0,X_1,X_2,X_3)$ est un polyn\^ {o}me homog\` ene de degr\'e $4$ dans $\mathbb P^3$, sans racines multiples et si $X$ d\'esigne le lieu d'annulation de $F$, alors $X$ est une surface $K3$ et la $2$-forme holomorphe est le r\'esidu de
$$\frac {dx_1\wedge dx_2 \wedge dx_3}{F(x_1,x_2,x_3)}$$
 o\` {u} $x_i:=\frac {X_i}{X_0}$.

 Le groupe $H_2(X,\mathbb {Z})$ est
libre de rang $22$ et l'accouplement d'intersection ou cup-produit munit 
$H_2(X,\mathbb {Z})$ d'une forme bilin\'eaire sym\'etrique unimodulaire, paire, de rang 
$22$, de signature $(3,19)$ telle que 

$$
H_2(X,\mathbb {Z})\simeq U_2^3\bot (-E_8)^2:=\mathcal {L} $$

o\`{u} $U_2$ est le r\'eseau hyperbolique de rang $2$ et $E_8$ le r\'eseau
unimodulaire d\'efini positif de rang $8$. Le r\'eseau 
$\mathcal {L}$
est appel\'e le r\'eseau $K3$. 

 Le groupe de Picard de $X$, not\'e $Pic X$, form\'e des diviseurs de 
$X$ modulo

 l'\'equivalence lin\'eaire, v\'erifie
$$
PicX\subset H^2(X,\mathbb {Z})\simeq Hom(H_2(X,\mathbb
{Z}),\mathbb {Z})$$

et $Pic X$ est param\'etr\'e par les cycles alg\'ebriques. C'est un groupe
ab\'elien libre de type fini, sans torsion, d'o\`{u} 

$$
Pic X \simeq \mathbb {Z} ^{\rho(X)}.$$

L'entier $\rho (X)$, appel\'e nombre de Picard de $X$, v\'erifie 
$$
1\leq \rho (X) \leq 20.$$

Le groupe $T(X):=(Pic X)^{\bot}$ des cycles transcendants a une structure de
r\'eseau de dimension $22-\rho (X)$. Il est appel\'e le r\'eseau
transcendant. 

 Si $X$ d\'esigne  une surface $K3$, $\mathcal {L}$ son r\'eseau $K3$ et $\alpha$ l'isomorphisme
$$\alpha:H_2(X,\mathbb {Z})\rightarrow \mathcal{L}, $$
le couple $(X,\alpha)$ est appel\'e surface $K3$ "marqu\'ee". 

 Si $\{ \gamma_1,...,\gamma_{22} \}$ d\'esigne une $\mathbb {Z}$-base de $H_2(X,\mathbb {Z})$ et $\omega$ une 2-forme holomorphe sur $X$, l'int\'egrale 
$\int_{\gamma _i}\omega$ est appel\'ee une p\'eriode de $X$ et v\'erifie 
$\int_{\gamma}\omega=0$ pour tout $\gamma \in Pic X$. 

 Si 
$\mathcal {M}$
est un sous-r\'eseau primitif de 
$\mathcal {L}$
( i. e. $\mathcal {L}/\mathcal {M}$ libre) de rang $1+t$, de signature $(1,t)$, le couple
$(X_{\mathcal {M}}
,\phi_\alpha)$ o\`{u} $X_\mathcal M$ est une surface $K3$ alg\'ebrique et $
\phi_\alpha=\alpha ^{-1}_{\mid \mathcal {M}}:\mathcal {M}\rightarrow Pic X_\mathcal M$ est une
isom\'etrie de r\'eseaux, est appel\'e surface $K3$, 
$\mathcal {M}$
-polaris\'ee. 

 On peut montrer l'existence de l'espace des modules des surfaces $K3$, 
$\mathcal {M}$
-polaris\'ees et pseudo-amples (i. e. dont le plongement $\phi_\alpha$
contient une classe de diviseurs pseudo-amples). Cet espace de modules est
ind\'ependant du marquage; il est not\'e $\mathscr{M}_{K3,\mathcal{M}}$. 

Supposons d\'esormais $\mathcal{M} \subset \mathcal{L}$ avec $
\mathcal{M}$ de rang $19$. 

 Si $\mathcal{M}\simeq U_2 \bot (-E_8)^2\bot \langle -2 \rangle>$, par un
th\'eor\`eme de Dolgachev,

 on a l'isomorphisme 
$$\mathscr{M}_{K3,\mathcal{M}}\simeq \mathcal{H}/\Gamma_0(N)^*$$
o\`{u} $\mathcal{H}$ d\'esigne le demi-plan de Poincar\'e, 

 $$
\Gamma_0(N)=\{ \left( 
\begin{matrix}
a & b \\ 
c & d
\end{matrix}
\right) \in Sl_2(\mathbb{Z}) / c\equiv 0 (N)\} 
$$

 et 
$$
\Gamma_0(N)^*=\Gamma_0(N)^{+}w_N
$$

o\`{u} $w_N$ d\'esigne l'involution de Fricke 
\begin{equation*}
w_N= \left( 
\begin{matrix}
0 & -\frac {1}{\sqrt{N}} \\ 
\sqrt{N} & 0
\end{matrix}
\right) 
\end{equation*}
Le groupe $\Gamma_0(N)^*$ est en outre de genre $0$. 

 Le groupe $\Gamma_0(N)^*$ (ou ses sous-groupes d'indice fini) peut
s'identifier au groupe de monodromie de l'\'equation diff\'erentielle de
Picard-Fuchs d'un pinceau de surfaces $K3$, $\mathcal {M}$-polaris\'ees (
pour la d\'efinition de l'\'equation diff\'erentielle de Picard-Fuchs, voir
ci-dessous). 

 Rappelons que $\mathcal{H}/\Gamma_0(N)^*$ est l'espace des modules
des couples $(E,C_N)$ des courbes elliptiques isog\`enes, \`a groupe
d'isog\'enie cyclique $C_N$, modulo l'involution de Fricke 
$w_N((E,C_N))=(E/C_N,EN
)$.Or le th\'eor\`eme de Dolgachev prouve que $\mathcal{H}/\Gamma_0(N)^*$
est \'egalement l'espace des modules des surfaces $K3$, 
$\mathcal {M}$
-polaris\'ees. On comprend donc qu'il puisse exister une relation entre les
surfaces $K3$, $\mathcal{M}$-polaris\'ees de nombre de Picard $19$ et les
courbes elliptiques. C'est pr\'ecis\'ement ce que met en \'evidence un
th\'eor\`eme de Morrison qui montre qu'une surface $K3$,
$\mathcal {M}$-polaris\'ee, de nombre de Picard $19$ poss\`ede
une structure de Shioda-Inose, i. e. il existe une surface
ab\'elienne $A:=E\times E/C_N$, une surface de Kummer
$Y=Kum(A/\pm 1)$ et une involution canonique $\iota$ sur $X$
telle que $X/\langle \iota \rangle $ soit birationnellement isomorphe \`a $Y$.

Consid\'erons maintenant une famille \`a $1$ param\` etre $X_z$ de
surfaces $K3$ param\'etr\'ee par $B:=\mathbb {P}^1 \backslash \{z / X_z \,\,{\hbox {singuli\` ere}} \}$ et soit $\omega _z$ l'unique
2-forme diff\'erentielle holomorphe sur $X_z$ (unique \`a un
scalaire pr\`es). Soit $z_0 \in B$ et $\pi (B,z_0)$ le groupe
fondamental. L'image de la repr\'esentation de monodromie
$$\pi (B,z_0)\rightarrow Aut(\mathbb {P} (H_2(X_z,\mathbb {Z}))$$
est le groupe de monodromie $G$ de la famille $\{X_z\}_{z\in B}$.
On d\'efinit \'egalement l'application de p\'eriode
$$\begin{aligned}
B & \rightarrow &\mathbb {P} ^{21} /G \\
z & \mapsto & \left[\int_{\gamma _{1 z}} \omega_z :... :\int_{\gamma_{22 z}} \omega_z \right] \end{aligned}
$$
On montre alors que si $X_z$ est une famille \`a un param\`etre
de surfaces $K3$, de nombre de Picard g\'en\'erique $r$, alors
les p\'eriodes de $X_z$ satisfont une \'equation diff\'erentielle
de Picard-Fuchs d'ordre $k=22-r$.

Dans nos exemples, nous aurons $k=3$.

\section{Preuve des th\'eor\`emes}
\subsection{Preuve du th\'eor\`eme $1.1$}

1) Rappelons d'abord les r\'esultats de Peters et Stienstra \cite{PS}
sur la famille $X_k$ de surfaces $K3$ dont une \'equation affine est 
$$ x+\frac
{1}{x}+y+\frac{1}{y}+z+\frac{1}{z}-k=0.$$
Une telle famille $\{X_k\}_k$, $k\in \mathbb {P}^1 \backslash
\{\infty ,\pm 2, \pm 6\}$ a un nombre de Picard g\'en\'erique
$19$, est $\mathcal {M}_k$-polaris\'ee avec
$$\mathcal{M}_k\simeq U_2 \bot (-E_8)^2\bot \langle -12 \rangle.$$
Son r\'eseau transcendant v\'erifie $T_k \simeq U_2 \bot \langle 12 \rangle$ et
l'\'equation diff\'erentielle de Picard-Fuchs associ\'ee \`a la
famille est 
$$(k^2-4)(k^2-36)y'''+6k(k^2-20)y''+(7k^2-48)y'+ky=0.$$
Si l'on pose
$$t(\tau)=\left ( \frac {\eta (\tau) \eta (6 \tau)}{\eta (2\tau)
\eta (3 \tau)}\right )^6=e^{\pi i \tau} \prod
_
{n=1\,(n,6)=1}^{\infty} {(1-e^{2\pi i \tau n})^6}$$

o\`{u} $\tau \in \mathcal {H}$, on peut montrer que
$$t(\frac {a\tau +b}{c\tau +d})=t(\tau) \, \, \, \, \forall  \left( 
\begin{matrix}
a & b \\ 
c & d
\end{matrix}
\right) \in \Gamma_1(6,2)^* \subset \Gamma _0(12)^* +12$$
o\`{u}
$$\Gamma _1 (6)= \{  \left( 
\begin{matrix}
a & b \\ 
c & d
\end{matrix}
\right) \in Sl_2(\mathbb {Z}) \,\,\, / \, \, \,a\equiv d \equiv 1
\, \, (6) \, \, \, c \equiv 0 \,\, (6) \} $$
$$\Gamma_1(6,2)=\{ \left( 
\begin{matrix}
a & b \\ 
c & d
\end{matrix}
\right) \in \Gamma_1(6) \,\,\, c\equiv 6b \,\,\,(12)\}$$
et $\Gamma_1 (6,2) ^*$ est le groupe engendr\'e par
$\Gamma_1(6,2)$ et l'involution de Fricke $w_6$.

En outre, $t$ est un Hauptmodul pour $\Gamma _1  (6,2)^*$, i. e.
induit un isomorphisme entre $\mathcal {H}^*=\mathcal {H} \cup
\mathbb {Q} \cup \{i \infty\}/\Gamma_1(6,2)^*$ et $\mathbb
{P}_1$.

On peut montrer que pour $\tau =i \infty$, on a $t=0$, pour $\tau
= \pm 1/2$, on a $t=\infty$, pour $\tau =i/\sqrt 6$, on a
$t=3-2\sqrt 2$ et pour $\tau=\pm 2/5 +i/5\sqrt 6$, on a
$t=3+2\sqrt 2$.

De plus, si $k=t+\frac {1}{t}$, l'\'equation de Picard-Fuchs en
la variable $t$ poss\`ede une base de solutions de la forme
$G(\tau), \tau G(\tau), \tau^2G(\tau)$ avec $G(\tau)=\eta(\tau) \eta(2\tau)\eta(3\tau)\eta(6\tau)$.

On a \'egalement
$$G(t)=\sum_{n \geq 0}
v_nt^{2n+1}\,\,\,\,\,\ \,\,\,\,\,\,\vert t \vert <3-2\sqrt2$$
avec
$$v_n=\sum_{k=0}^n \left ( \begin{matrix}
 n\\k \end{matrix}
\right )^2 \left (
\begin{matrix} n+k\\k \end{matrix}
\right)^2.$$ 
Nous allons maintenant prouver la premi\`ere assertion du
th\'eor\`eme 1.

Par d\'efinition,
$$m(P_k)=\frac {1}{(2\pi i)^3} \int_{\mathbb {T}^3} log \vert
k-(x+\frac {1}{x}+y+\frac {1}{y}+z+\frac {1}{z})\vert \frac
{dx}{x} \frac {dy}{y} \frac {dz}{z}.$$
Et pour $k>6$,
$$ \begin{aligned}
\frac {dm(P_k)}{dk} &=\frac {1}{(2\pi i)^3}\frac {1}{k} \int_{\mathbb
{T}^3} \frac {1}{1-\frac {1}{k}(x+\frac {1}{x}+y+\frac {1}{y}+z+\frac {1}{z})} \frac
{dx}{x} \frac {dy}{y} \frac {dz}{z} \\
                    &=\sum_{m\geq 0} a_m \frac {1}{k^{2m+1}} \end{aligned}$$
avec
$$a_m=\sum _{p+q+r=m} \frac {(2m)!}{(p!q!r!)^2}.$$
Or $\frac {dm(P_k)}{dk}$ est une p\'eriode et donc v\'erifie
l'\'equation diff\'erentielle de Picard-Fuchs. En faisant alors
le changement de variable $k=t+\frac {1}{t}$, on obtient
$$\frac {dm(P_k)}{dk}=\sum _{n\geq 0}v_nt^{2n+1}=G(\tau)$$
soit
$$dm(P_k)=-G(\tau)\frac {dt}{t} \frac {1-t^2}{t}.$$
Comme $ -G(\tau) q \frac {dt}{t} \frac {1-t^2}{t}$ est une forme
modulaire de poids $4$ pour $\Gamma_1(6,2)^*$,
$$\begin{aligned}
-G(\tau) q \frac {dt}{t} \frac {1-t^2}{t} &=-\frac {1}{2} +4q+20q^2+148q^3+148q^4+504q^5\\
                                           &+740q^6+1376q^7+1172q^8+0(q^8), \end{aligned}$$
 on va chercher \`a l'\'ecrire sous la forme 
$$\alpha E_4(\tau)+\beta E_4(2\tau)+\gamma E_4(3\tau)+\delta
E_4(6\tau)$$
o\`{u}
$$\begin{aligned}
 E_4(\tau) &=1+240\sum_{n\geq 1}(\sum_{d\mid n}d^3)q^n\\
           &=1+240(q+9q^2+28q^3+73q^4+126q^5+252q^6+344q^7\\
           &+585q^8+757q^9+1134q^{10}+...)
\end{aligned}$$
En calculant alors suffisamment de termes dans leur
$q$-d\'eveloppement, on voit que ces deux formes co\"\i ncident
pour 
$$\alpha = \frac {4}{240} \,\,\,\,\beta =- \frac
{16}{240}\,\,\,\,\gamma=\frac {36}{240}\,\,\,\,\delta=-\frac
{144}{240}.$$
On a donc
$$ \begin{aligned}
dm(P_k) &=(\frac {4}{240} E_4(q)-\frac {16}{240} E_4(q^2)+\frac
{36}{240}E_4(q^3)-\frac {144}{240}E_4(q^6))\frac {dq}{q}\\
        &=(-\frac {1}{2q} +4q+20q^2+...)dq
\end{aligned}$$
En int\'egrant entre $k$ et l'infini, on trouve alors
$$m(P_k)=\Re (-\pi i \tau +\sum_{n\geq 1} (\sum_{d\mid n} d^3)(4\frac {q^n}{n}-8\frac {q^{2n}}{n}+12\frac {q^{3n}}{n}-24\frac {q^{6n}}{n})).$$ 
2) Rappelons maintenant les r\'esultats de Verrill \cite{V}
sur la famille $Y_k$ de surfaces dont une \'equation affine est 
$$  (1+x+xy+xyz)(1+z+zy+zyx)-(k+4)xyz=0.$$
Une telle famille $\{Y_k\}_k$, $k\in \mathbb {P}^1 \backslash
\{\infty ,-4,12,0\}$ a un nombre de Picard g\'en\'erique
$19$, est $\mathcal {M}_k$-polaris\'ee avec
$$\mathcal{M}_k\simeq U_2 \bot (-E_8)^2\bot <-6>.$$
Son r\'eseau transcendant v\'erifie $T_k \simeq U_2 \bot <6>$ et
l'\'equation diff\'erentielle de Picard-Fuchs associ\'ee \`a la
famille est \cite{VY}
$$k(k+4)(k-12)y'''+6(k^2-7k-12)y''+\frac {7k^2-12k-96}{k+4}
y'+\frac {k}{k+4} y=0.$$
Si l'on pose
$$t(\tau)= \frac {\eta (3\tau)^4 \eta (12 \tau)^8 \eta (2\tau)^{12}}{\eta (\tau)^4
\eta (4 \tau)^8\eta (6\tau)^{12}}$$

o\`{u} $\tau \in \mathcal {H}$, on peut montrer que
$$t(\frac {a\tau +b}{c\tau +d})=t(\tau) \, \, \, \, \forall  \left( 
\begin{matrix}
a & b \\ 
c & d
\end{matrix}
\right) \in  \Gamma _0(12) +12$$
o\`{u}
$\Gamma_0 (12)+12 $ est le groupe engendr\'e par
$\Gamma_0(12)$ et l'involution de Fricke $w_{12}$.

En outre, $t$ est un Hauptmodul pour ce groupe.

On peut montrer que pour $\tau =i \infty$, on a $t=0$ et pour $\tau=i0$, on a
$t=-1$.

De plus, si $k=-(t+\frac {1}{t}+2)$, l'\'equation de Picard-Fuchs en
la variable $t$ poss\`ede une base de solutions de la forme
$G(\tau), \tau G(\tau), \tau^2G(\tau)$ avec $G(\tau)=\frac {\eta(2\tau)^4 \eta(6\tau)^4}{\eta(\tau)^2\eta(3\tau)^2}$.

On a \'egalement
$$G(t)=\sum_{n \geq 1}
v_nt^n$$
avec
$$v_n=\sum_{m=0}^{n-1}\sum_{p+q+r+s=m}(-1)^m \left ( \begin{matrix}
 n+m\\2m+1 \end{matrix}
\right ) \left (
\begin{matrix} m!\\p!q!r!s! \end{matrix}
\right)^2.$$ 
Nous allons maintenant prouver la deuxi\`eme assertion du
th\'eor\`eme 1.1.

Comme pr\'ec\'edemment on trouve
$$dm(k)=-G(\tau)\frac {dt}{t}\frac {1-t^2}{t}$$
et 
$$-G(\tau)q\frac {dt}{t}\frac {1-t^2}{t}=-\frac
{2}{240}E_4(\tau)+\frac {32}{240}E_4(2\tau)+\frac
{18}{240}E_4(3\tau)-\frac {288}{240}E_4(6\tau)$$
car$-G(\tau)\frac {dt}{t}\frac {1-t^2}{t}$ est une forme
modulaire de poids $4$ pour $\Gamma _0(12)+12$.

En int\'egrant entre $k$ et l'infini, on trouve alors le
r\'esultat annonc\'e.

\subsection{Preuve du th\'eor\`eme $1.2$}

Les \'etapes de la d\'emonstration sont semblables \`a celles
d\'evelopp\'ees dans \cite{Ber}.

1) Partant de la relation
$$
m(P_k)=\Re \{-\pi i \tau+\sum_{n\geq 1}(\sum_{d\mid n}d^3)(\frac
{4q^n}{n}-\frac {8q^{2n}}{n}+\frac {12q^{3n}}{n}-\frac {24q^{6n}}{n})\},
$$
on pose $n=dn'$, puis gr\^{a}ce \`a la relation
$$D^2(Li _3(q^{jd}))=j^2d^2Li_1(q^{jd})
\,\,\,\,\,\,\,j=1,2,3,6\,\,\,\,\,\,D=q\frac {d}{dq},$$
on obtient
$$m(k)=\Re \{ -\pi i\tau +4D^2(\sum_{d\geq 1} Li_3(q^d)-\frac
{1}{2} Li_3(q^{2d})+\frac {1}{3} Li_3(q^{3d})-\frac
{1}{6}Li_3(q^{6d})\}.$$
Notons alors
$$L_j(x)=\sum_{d\geq 1} Li_3(q^{jd}x)$$
et
$$H_j(x)=L_j(x)+L_j(\frac {1}{x}).$$
On peut montrer qu'il existe des constantes complexes non nulles
$A$, $B$, $C$, $D$ telles que
$$\begin{aligned}
K(x) &=H_1(x)-\frac {1}{2}H_2(x)+\frac {1}{3}H_3(x)-\frac
{1}{6}H_6(x)+A\log (x)^4+B\log (x)^3\\
     &+C\log (x)^2+D\log (x) \end{aligned}$$

soit invariant par la transformation
$$x\mapsto q^6x.$$
En effet 
$$H(x)=H_1(x)-\frac {1}{2}H_2(x)+\frac {1}{3}H_3(x)-\frac
{1}{6}H_6(x)$$
se transforme en
$$\begin{aligned} 
H(x) &-(Li_3(qx)-Li_3(\frac {1}{qx}))-\frac {1}{2}(Li_3(q^2x)-Li_3(\frac
{1}{q^2x}))\\
     &-\frac {4}{3}(Li_3(q^3x)-Li_3(\frac {1}{q^3x}))-\frac {1}{2}(Li_3(q^4x)-Li_3(\frac
{1}{q^4x}))\\
     &-(Li_3(q^5x)-Li_3(\frac
{1}{q^5x}))-\frac {2}{3}Li_3(q^6x)+\frac {2}{3} Li_3(\frac {1}{x})
\end{aligned}$$
et l'on a la formule
$$Li_3(z)-Li_3(\frac {1}{z})=-\frac {(2i\pi)^3}{6} B_3(\frac
{\log z}{2i\pi})$$
o\`{u} $B_3$ d\'esigne le polyn\^{o}me de Bernouilli
$$B_3(X)=X^3-\frac {3}{2}X^2+\frac {1}{2}X.$$
Par suite $H(x)$ se transforme en 
$$H(x)+A'\log (x)^3+B'\log (x)^2+C'\log (x)+D'$$
et $K(x)$ en
$$\begin{aligned} 
K(x) & +(4\lambda A+A')\log (x)^3+(6\lambda ^2A+3\lambda B+B')\log
(x)^2\\
     & +(4A\lambda ^3+3B\lambda ^2+2C\lambda +C')\log (x)+A\lambda
^4+B\lambda ^3+C\lambda ^2+D\lambda +D',
\end{aligned} $$
avec $\lambda=6\log q$.

Il suffit alors de d\'eterminer $A$, $B$, $C$, $D$ en fonction de
$A'$, $B'$, $C'$, $D'$.

Comme 
$$m(k)=\Re (-\pi i \tau +2 D^2(K(1)))$$ et comme $K(e^{2\pi i\xi \tau})$ est
invariant par le changement $\xi \mapsto \xi +6$ d'apr\`es ce qui
pr\'ec\`ede, on va d\'evelopper $K(e^{2\pi i \xi \tau})$ en
s\'erie de Fourier. Le d\'eveloppement en s\'erie de $K(1)$ sera
obtenu pour $\xi =0$.

Expliquons le calcul du d\'eveloppement sur $L_1(x)$.
On va \'ecrire
$$\begin{aligned}
L_1(e^{2\pi i \xi \tau})&=\sum_{d\geq 1}\sum_{m\geq 1}\frac
{e^{2\pi i \tau m(d+\xi)}}{m^3}\\
                        &=\sum _{m\geq 1} \frac {1}{m^3} (\sum
_{d\equiv 1\,\,(6)}e^{2\pi i \tau m(d+\xi)}+...+
\sum_{d\equiv 6\,\,(6)}e^{2\pi i
\tau m(d+\xi)})
\end{aligned} $$
Ensuite on calcule
$$I_{n,h}=\sum_{k\geq 0}\frac {1}{6}\int_{-h}^{6-h}e^{2\pi i
m\tau(6k+h+\xi)}e^{-2\pi in\frac {\xi}{6}}d\xi$$
Posant alors $\xi'=6k+h+\xi$, il vient
$$I_{n,h}=\frac {1}{6}e^{\frac {2\pi i
nh}{6}}\int_0^{+\infty}e^{2\pi i\xi'(m\tau-\frac {n}{6})}d\xi'.$$
soit
$$I_{n,h}=\frac {-1}{6}e^{\frac {2\pi i
nh}{6}}\frac {1}{2\pi i (m\tau -\frac {n}{6})}.$$
On en d\'eduit alors
$$\frac {1}{6}\int_{\hbox{p\'eriode}}L_1(e^{2\pi i\tau
\xi})e^{-2\pi in\frac {\xi}{6}}d\xi=-\frac {1}{2\pi i}\sum_{n\geq
1}\frac {1}{m^3}\frac {1}{(m\tau-\kappa)}\,\,\,\,{\hbox {si}}\,\,\,\,n=6\kappa$$
et $0$ sinon.

D'o\`{u}
$$\begin{aligned}
K(1)=-\frac {1}{2\pi i}\sum_{\kappa}(& \sum_{m\geq 1}\frac {1}{m^3}(\frac
{1}{m\tau-\kappa}+\frac {1}{m\tau+\kappa})\\
                              &-\frac {1}{2}\sum_{m\geq 1}\frac {1}{m^3}(\frac
{1}{2m\tau-\kappa}+\frac {1}{2m\tau+\kappa})\\
                              &\frac {1}{3}\sum_{m\geq 1}\frac {1}{m^3}(\frac
{1}{3m\tau-\kappa}+\frac {1}{3m\tau+\kappa})\\
                              &-\frac {1}{6}\sum_{m\geq 1}\frac {1}{m^3}(\frac
{1}{6m\tau-\kappa}+\frac {1}{6m\tau+\kappa}))
\end{aligned}$$ 
Finalement

$$\begin{aligned}
m(k) &=\Re(-\pi i \tau-\frac {4i}{8\pi^3}\sum_{\kappa,m\neq 0}\frac
{1}{m}(\frac {1}{(m\tau+\kappa)^3}-2\frac {1}{(2m\tau+\kappa)^3}+3\frac {1}{(3m\tau+\kappa)^3}\\
   &-6\frac {1}{(6m\tau+\kappa)^3})) \end{aligned}$$

D'o\`{u} le r\'esultat annonc\'e, puisque
$$\Im \frac {1}{(m\tau+\kappa)^3}=-m\Im \tau (2\Re(\frac
{1}{(m\tau+\kappa)^3(m\bar{\tau}+\kappa)})+\frac
{1}{(m\tau+\kappa)^2(m\bar {\tau}+\kappa)^2})$$
et
$$\frac {\Im \tau }{8\pi^3}6\times 120 \sum_{k\geq 1}\frac
{1}{k^4}=\pi \Im \tau.$$

\section{Quelques applications}
\subsection{Valeur approch\'ee de la mesure}
La formule de Jensen permet d'exprimer la mesure de Mahler d'un polyn\^{o}me de trois variables \` a l'aide d'une int\'egrale double. Cependant les m\'ethodes d'int\'egration num\'erique demanderaient beaucoup de temps pour obtenir une pr\'ecision de $10^{-18}$ par exemple. Les formules pr\'ec\'edentes expriment la mesure de Mahler \` a l'aide de s\'eries rapidement convergentes.

Boyd et Mossinghoff ont ainsi trouv\'e la valeur approch\'ee de $P_1$

$$M(x+y+z+1+\frac {1}{x}+\frac {1}{y}+\frac {1}{z})=1,4483035845491699038...$$

J'ai de m\^ {e}me calcul\'e une valeur approch\'ee de la mesure du polyn\^ {o}me $Q_1$
$$M(Q_1)=1,435170000343077634...$$

Ces deux mesures sont parmi les plus petites mesures connues pour les polyn\^ {o}mes de $3$ variables dont la mesure ne se r\'eduit pas \` a celle d'un polyn\^ {o}me de deux variables comme par exemple
$$\begin{aligned}  
M(x+y+z+\frac {1}{x}+\frac {1}{y}+\frac {1}{z})&=M(1+x+y)\\
                                               &=1,3813564445184977933...
\end{aligned}
$$
\subsection{Mesure et s\'erie $L$ de Hecke}
Pour certaines valeurs de $k$ correspondant \` a des $\tau$ quadratiques donc \` a des surfaces $K3$ singuli\` eres (i.e. de nombre de Picard $20$) ayant une structure de Shioda-Inose li\'ee \` a une courbe elliptique \` a multiplication complexe par un ordre d'un corps quadratique imaginaire, la mesure de Mahler s'exprime \` a l'aide d'une s\'erie $L$ de Hecke d'un ordre. Ceci se produit par exemple pour $P_k$ avec $k=2,3,6$ et $Q_k$ avec $k=-4,12$.

Soit $K=\mathbb Q(\sqrt{d})$ un corps quadratique imaginaire d'anneau des entiers $O_K$ et de discriminant $D$. Un Gr\"{o}ssencharacter $\phi$ de poids $k\geq 2$, de conducteur $\Lambda$, o\` {u} $\Lambda$ est un id\'eal de $O_K$ est ainsi d\'efini. Un homomorphisme $\phi :I(\Lambda) \rightarrow \mathbb C^\times $ satisfaisant
$$\phi (\alpha O_K)=\alpha^{k-1} \,\,\,\,\,_{\hbox {pour} }\,\,\,\,\,\alpha \equiv 1 \,\,\,\,\hbox{mod}\,\,\,\Lambda$$
est appel\'e Gr\"{o}ssencharacter de Hecke de poids $k$ et conducteur $\Lambda$. La s\'erie $L$ de Hecke induite par le Gr\"{o}ssencharacter de Hecke est d\'efinie par
$$L(\phi,s):=\sum_P \frac {\phi (P)}{N(P)^s} =\sum_{n=1}^{\infty} \frac {a(n)}{n^s}$$

o\` {u} $ N(P)$ est la norme de l'id\'eal $P$ et la somme prise surles id\'eaux premiers 
$P\subset O_K$ premiers \` a $\Lambda$. 
En rempla\c {c}ant $O_K$ par un ordre $R$ du corps quadratique, on d\'efinirait de m\^{e}me la s\'erie $L$ de Hecke d'un ordre pour un Gr\"{o}ssencharacter $\phi$ attach\'e \` a l'ordre. Dans tous les cas on pr\'ecisera $L_{\mathbb Q(\sqrt{d})}$ ou $L_R$.

\begin{theo}
Les Gr\"{o}ssencharacter  de ce th\'eor\` eme sont tous de poids $3$.
$$m(P_0)=d_3:=\frac {3\sqrt {3}}{4\pi} L(\chi_{-3},2).$$

$$m(P_2)=\frac {16\sqrt{2}}{\pi^3}L_{\mathbb Q(\sqrt{-2})}(\phi ,3)$$

$$m(P_3)=\frac {15\sqrt{15}}{2\pi^3}L_{\mathbb Q(\sqrt{-15})}(\phi ,3)$$
o\` {u} $\phi(P)=-\omega$ si $P=(2,\omega)$ et $\omega=\frac {1+\sqrt{-15}}{2}$, $P$ d\'esignant un repr\'esentant de la deuxi\` eme classe d'id\'eaux du corps de nombres $\mathbb Q(\sqrt{-15})$ de nombre de classe $2$.

$$m(P_6)=\frac {24\sqrt{6}}{\pi^3}L_{\mathbb Q(\sqrt{-6})}(\phi ,3)$$
o\` {u} $\phi(P)=-2$ si $P=(2,\sqrt{-6})$, $P$ d\'esignant un repr\'esentant de la deuxi\` eme classe d'id\'eaux du corps de nombres $\mathbb Q(\sqrt{-6})$ de nombre de classe $2$.

$$m(Q_0)=\frac {12\sqrt{3}}{\pi^3}L_R(\phi ,3)$$
pour l'ordre $(1,2\sqrt{-3})$ de nombre de classe $1$.

$$m(Q_{12})=4m(Q_0).$$

\end{theo}

{\it Preuve}
La preuve utilise le th\'eor\`eme 2 et la formule

$$L_F(\phi,s)=\sum_{cl(P)} \frac {\phi(P)}{N(P)^{2-s}}Z(2,P,s)$$

o\` {u}

$$Z(2,P,s)=\frac {1}{2}\sum'_{\lambda \in P} \frac {\bar {\lambda}^2}{(\lambda \bar {\lambda})^s}$$

est la s\'erie de Hecke partielle.

Pour cela, \'ecrivons $m(P_k)$ sous une autre forme.

Posons 
$$D_{j\tau}=(mj\tau+\kappa)(mj\bar {\tau}+\kappa).$$
Alors
$$\begin{aligned}
m(P_k)=\frac {\Im \tau}{8\pi^3}\sum'_{m,\kappa} & [-4\frac
{(m(\tau+\bar {\tau})+2\kappa)^2}{D_{\tau}^3}+\frac {4}{D_{\tau}^2}\\
                                             & +16\frac
{(2m(\tau+\bar {\tau})+2\kappa)^2}{D_{2\tau}^3}-\frac
{16}{D_{2\tau}^2}\\
                                             & -36\frac
{(3m(\tau+\bar {\tau})+2\kappa)^2}{D_{3\tau}^3}+\frac
{36}{D_{3\tau}^2}\\
                                             & +144\frac
{(6m(\tau+\bar {\tau})+2\kappa)^2}{D_{6\tau}^3}-\frac
{144}{D_{6\tau}^2}]
\end{aligned}$$ 
Si $k=6$, on a $\tau=\frac {i}{\sqrt {6}}$ et 
$$D_\tau=\frac {1}{6}(m^2+6\kappa^2)$$
$$D_{2\tau}=\frac {1}{3}(2m^2+3\kappa^2)$$
$$D_{3\tau}=\frac {1}{2}(3m^2+2\kappa^2)$$
$$D_{6\tau}=(6m^2+\kappa^2).$$
D'o\`{u}
$$m(P_6)=\frac {24\sqrt {6}}{\pi^3}[\frac
{1}{2}\sum'_{m,\kappa}(\frac
{m^2-6\kappa^2}{(m^2+6\kappa^2)^3}+\frac
{3\kappa^2-2m^2}{(3\kappa^2+2m^2)^3})].$$
Or dans le corps ${\mathbb Q}(\sqrt {-6})$, de discriminant $-24$,
il y a $2$ classes d'id\'eaux, celle $\mathcal A_0=\{(\lambda)\}$
des id\'eaux principaux et la classe $\mathcal A_1=\{(\lambda
)\mathcal P\}$ o\`{u} ${\mathcal P}=(2,\sqrt {-6})$.

Si l'on d\'efinit le caract\`ere de Hecke par
$$\psi ((\lambda))=\lambda ^2$$
pour $\lambda =m+\sqrt {-6}\kappa$
et $$\psi (\mathcal P)=-2,$$
on obtient la formule annonc\'ee.

\bigskip

Si $k=2$, on a $\tau=-\frac {1}{3}+i\frac {\sqrt {2}}{6}$ et 
$$D_\tau=\frac {1}{6}((m-2\kappa)^2+2\kappa^2)$$
$$D_{2\tau}=\frac {1}{3}(2(m-\kappa)^2+\kappa^2)$$
$$D_{3\tau}=\frac {1}{2}(2(\kappa-m)^2+m^2)$$
$$D_{6\tau}=(\kappa -2m)^2+2m^2.$$
D'o\`{u}
$$m(P_2)=\frac {4\sqrt {2}}{\pi^3}\sum'_{m,\kappa}(\frac
{m^2-6\kappa^2}{((m-2\kappa)^2+2\kappa^2)^3}+\frac
{3\kappa^2-2m^2}{(2(m-\kappa)^2+\kappa^2)^3}).$$
En posant $m-2\kappa=l$ dans la premi\`ere fraction puis
$m-\kappa=\kappa'$ et $\kappa=l'$ dans la seconde, on trouve
$$m(P_2)=\frac {16\sqrt {2}}{\pi^3}\frac
{1}{2}\sum'_{\kappa,l}\frac {l^2-2\kappa^2}{(l^2+2\kappa^2)^3}$$
c'est-\`a-dire la formule annonc\'ee.

\bigskip

Pour $k=3$, on a $\tau=\frac {-3+\sqrt {-15}}{12}$.
$$D_\tau=\frac {1}{6}(m^2-3\kappa m+6\kappa^2)$$
$$D_{2\tau}=\frac {1}{3}(3\kappa^2-3\kappa m+2m^2)$$
$$D_{3\tau}=\frac {1}{2}(2\kappa^2-\kappa m+3m^2)$$
$$D_{6\tau}\kappa^2-3\kappa m+6m^2.$$
Or il existe exactement deux formes binaires quadratiques r\'eduites de discriminant $-15$, \` a savoir $(1,1,4)$ et $(2,1,2)$. A l'aide d'un changement de variables, on va exprimer $\Delta_{j\tau}$ en fonction de ces formes et en d\'eduire $m(P_3)$.

Pour $\Delta_{\tau}$ on pose $m=m'+2\kappa$, d'o\` {u} $D_{\tau}=\frac {1}{6}(m'^2+\kappa m'+4\kappa^2)$.

Pour $\Delta_{2\tau}$ on pose $m=m'+\kappa$, d'o\` {u} $D_{2\tau}=\frac {1}{3}(2m'^2+\kappa m'+2\kappa^2)$.

Pour $\Delta_{3\tau}$ on pose $\kappa=\kappa'+m$, d'o\` {u} $D_{3\tau}=\frac {1}{2}(2\kappa'^2+\kappa'm+2m^2)$.

Pour $\Delta_{6\tau}$ on pose $\kappa=\kappa'+2m$, d'o\` {u} $D_{6\tau}=\kappa'^2+\kappa'm+4m^2$.

D'o\` {u}, par abus de notation, puisque la sommation en $m'$ ou en $\kappa'$ est \'equivalente \` a la sommation en $m$ et $k$,

$$m(P_3)=\frac {15\sqrt{15}}{2\pi^3} \sum'_{m,\kappa} \frac {1}{2} (\frac {m^2+4m\kappa-2\kappa^2}{(m^2+\kappa m+4\kappa^2)^3}+\frac {m^2-4m\kappa-2\kappa^2}{(2m^2+\kappa m+2\kappa^2)^3}).$$

Comme $2m^2+\kappa m+2\kappa^2$ est sym\'etrique en $\kappa$ et $m$, on en d\'eduit

$$\sum'_{m,\kappa} \frac {m^2-4m\kappa-2\kappa^2}{(2m^2+\kappa m+2\kappa^2)^3}=\frac {1}{2}\sum'_{m,\kappa} \frac {-m^2-8m\kappa-\kappa^2}{(2m^2+\kappa m+2\kappa^2)^3}).$$

En posant $\kappa=\kappa'-m$, on va \'ecrire sous forme sym\'etrique
$$m^2+\kappa m+4\kappa^2=4\kappa'2-7m\kappa'+4m^2$$
et l'on obtient

$$\frac {m^2+4m\kappa-2\kappa^2} {(m^2+\kappa m+4\kappa^2)^3}=\frac {-5m^2+8m\kappa'-2\kappa'^2} {(4m^2-7\kappa'm+4\kappa'^2)^3}$$
puis

$$\begin{aligned}
\sum'_{m,\kappa'}\frac {-5m^2+8m\kappa'-2\kappa'^2} {(4m^2-7\kappa'm+4\kappa'^2)^3 } &=\frac {1}{2}\sum'_{m,\kappa'}\frac {-7m^2+16m\kappa'-7\kappa'^2}{(4m^2-7\kappa'm+4\kappa'^2)^3}\\
            &=\frac {1}{2}\sum'_{m,\kappa}\frac {2m^2+2m\kappa-7\kappa^2}{(m^2+\kappa m+4\kappa^2)^3}. \end{aligned}$$

Finalement,

$$m(P_3)=\frac {15\sqrt{15}}{2\pi^3} \sum'_{m,\kappa} \frac {1}{4} (\frac {2m^2+2m\kappa-7\kappa^2}{(m^2+\kappa m+4\kappa^2)^3}-\frac {m^2+8m\kappa+\kappa^2}{(2m^2+\kappa m+2\kappa^2)^3}).$$

Dans le corps $F=\mathbb Q(\sqrt{-15})$ il y a deux classes d'id\'eaux entiers, la classe des id\'eaux principaux dont un repr\'esentant est l'id\'eal $(1,\omega)=(1)$ et l'autre repr\'esent\'ee par l'id\'eal $P=(2,\omega)$. D\'efinissons le caract\` ere de Hecke de poids $3$ sur $P$ par $\phi (P)=-\omega$ puisque $P^2=(\omega)$. On obtient alors
$$Z(2,(1),s)=\frac {1}{2} \sum'_{\lambda \in (1)} \frac {\bar {\lambda}^2+\lambda^2}{(\lambda \bar{\lambda})^s}=\frac {1}{4} \sum'_{m,\kappa}\frac {2m^2+2m\kappa-7\kappa^2}{(m^2+m\kappa+4\kappa^2)^s}$$

car $\lambda=m+\kappa\omega$ et

$$\frac {1}{2}\frac {\phi(P)}{N(P)^{2-s}} \sum'_{\lambda \in P} \frac {\bar {\lambda}^2}{(\lambda \bar{\lambda})^s}=\frac {1}{2}\frac {\phi(\bar {P})}{N(\bar {P})^{2-s}} \sum'_{\bar {\lambda} \in \bar {P}} \frac {\lambda^2}{(\lambda \bar{\lambda})^s}=-\frac {1}{4}\frac {1}{2^{2-s}} \sum'_{\lambda \in P} \frac {\bar {\lambda}^2\omega+\lambda^2\bar {\omega}}{(\lambda \bar{\lambda})^s}$$

$$=-\frac {1}{4} \sum'_{m,\kappa} \frac {m^2+8m\kappa+\kappa^2}{(2m^2+\kappa m+2\kappa^2)^3},$$
puisque $\lambda =2m+\kappa\omega$.

Finalement
$$L_F (\phi,3)=\frac {1}{4} \sum'_{m,\kappa}(\frac {2m^2+2m\kappa-7\kappa^2}{(m^2+m\kappa+4\kappa^2)^3} -\frac {m^2+8m\kappa+\kappa^2}{(2m^2+\kappa m+2\kappa^2)^3}),$$

d'o\` {u} le r\'esultat annonc\'e.

\bigskip
Pour $k=0$, Boyd avait prouv\'e
$$m(P_0)=d_3:=\frac {3\sqrt {3}}{4\pi} L(\chi_{-3},2).$$

Nous allons retrouver ce r\'esultat autrement car la preuve met en \'evidence les relations entre fonction z\'eta d'un ordre et celle de l'ordre maximal.
Tout d'abord remarquons que $h(-12)=1$ et que la forme r\'eduite de discriminant $-12$ est $x^2+3y^2$. De m\^ {e}me $h(-3)=1$ et la forme r\'eduite de discriminant $-3$ est $x^2+xy+y^2$.

Pour $k=0$ on a $\tau=\frac {-3+\sqrt {-3}}{6}$ et

$$D_\tau=\frac {1}{3}(m^2-3\kappa m+3\kappa^2)$$
$$D_{2\tau}=\frac {1}{3}(3\kappa^2-6\kappa m+4m^2)$$
$$D_{3\tau}=\kappa^2-3\kappa m+3m^2$$
$$D_{6\tau}=\kappa^2-6\kappa m+12m^2.$$

Pour $\Delta_{\tau}$ on pose $m=m'+2\kappa$, d'o\` {u} $D_{\tau}=\frac {1}{3}(m'^2+\kappa m'+\kappa^2)$.

Pour $\Delta_{2\tau}$ on pose $\kappa=\kappa'+m$, d'o\` {u} $D_{2\tau}=\frac {1}{3}(m^2+3\kappa'^2)$.

Pour $\Delta_{3\tau}$ on pose $\kappa=\kappa'+2m$, d'o\` {u} $D_{3\tau}=m^2+\kappa' m+\kappa'^2$.

Pour $\Delta_{6\tau}$ on pose $\kappa=\kappa'+3m$, d'o\` {u} $D_{6\tau}=\kappa'^2+3m^2$.
Apr\` es simplification, on obtient

$$m(P_0)=\frac {3\sqrt {3}}{2\pi^3} \sum'_{m,\kappa}(\frac {4}{(m^2+3\kappa^2)^2}-\frac {1}{(m^2+m\kappa+\kappa^2)^2}).$$

Posons $K=\mathbb Q(\sqrt {-3})$ et d\'esignons par $R$ l'ordre de discriminant $-12$. Alors 
$m(P_0)$ s'\'ecrit
$$m(P_0)= \frac {3\sqrt {3}}{2\pi^3}(8\zeta_R(2)-6\zeta_K(2)).$$
Comme, par exemple \cite {Sh} 
$$8\zeta_R(2)=9\zeta_K(2)$$
et
$$\zeta_K(2)=\zeta(2)L(\chi_{-3},2),$$
on obtient bien 
$$m(P_0)=d_3.$$
\bigskip

Si $k=12$ dans la deuxi\` eme famille, alors $\tau=\frac {3+\sqrt {-3}}{6}$ et l'on obtient
$$m(Q_{12})=\frac {\sqrt {3}}{48\pi^3} \sum'_{\kappa,m} (8\times 4^2\times 3^2 \frac {\kappa^2-3m^2}{(\kappa^2+3m^2)^3}+18\frac {2m^2+2m\kappa-\kappa^2}{(m^2+\kappa m+\kappa^2)^3})$$

tandis que pour $k=0$ on a $\tau=\frac {3+\sqrt {-3}}{12}$ et 
$$m(Q_{0})=\frac {\sqrt {3}}{72\pi^3} \sum'_{\kappa,m} (8\times 4^3\times 3^2 \frac {\kappa^2-3m^2}{(\kappa^2+3m^2)^3}+18\times 16\frac {2m^2+2m\kappa-\kappa^2}{(m^2+\kappa m+\kappa^2)^3}).$$

Or d'apr\` es une remarque de Sebbar \cite {Se}
$$\sum'_{m,\kappa} \frac {2m^2+2m\kappa-\kappa^2}{(m^2+\kappa m+\kappa^2)^3}=0,$$

d'o\` {u}
$$m(Q_{12})=4m(Q_0)$$
relation conjectur\'ee par Boyd.\cite{Bo3}

En outre
$$m(Q_{0})=\frac {12\sqrt {3}}{\pi^3}\frac {1}{2} \sum'_{\kappa,m}  \frac {\kappa^2-3m^2}{(\kappa^2+3m^2)^3}=\frac {12\sqrt {3}}{\pi^3}L_R(\phi ,3)$$
o\` {u} $\phi$ d\'esigne le caract\` ere de Hecke de poids $3$ pour l'ordre $(1,2\sqrt{-3})$ d\'efini par $\phi((\alpha))=\alpha ^2$.
\qed

\bigskip

Boyd \cite {Bo3} a conjectur\'e une autre relation entre les mesures de Mahler de la deuxi\` eme famille:

$$(*)\,\,\,\,\,\,\,\,\,\, 2m(Q_{-36})=?4m(Q_{-6})+m(Q_0).$$
Nous allons montrer l'\'equivalence de cette relation avec une relation entre s\'eries $L$ de Hecke d'ordres d'un corps de nombres.
\begin{theo}
Soit $R=(1,2\sqrt{-3})$ et $R'=(1,\sqrt{-3})$ les deux ordres $R\subset R'$ du corps de nombres $\mathbb Q(\sqrt{-3})$ de discriminants respectifs $-48$ et $-12$, de nombre de classes respectif $2$ et $1$.

Soit $\phi_R$ (resp. $\phi_{R'}$) les caract\` eres de Hecke de poids $3$ d\'efinis par 
$$\phi_R(\alpha R)=\alpha^2      \,\,\,\,\,\,\,\,\,\phi_R (P)=-3 \,\,\,\,{\hbox {si}}\,\,\,\,\, P=(3,2\sqrt{-3})$$
$$\phi_{R'}(\beta R')=\beta^2.$$
La relation (*) entre mesures de Mahler est \'equivalente \` a la relation
$$L_{R'}(\phi_{R'},3)=L_R(\phi_R,3).$$

\end{theo}

\bigskip

{\it Preuve}
Pour $k=-6$, on a $\tau =\sqrt{-3}/6$ et
$$D_{\tau}=\frac {1}{12} (m^2+12\kappa^2)$$
$$D_{2\tau}=\frac {1}{3} (m^2+3\kappa^2)$$

$$D_{3\tau}=\frac {1}{4} (3m^2+4\kappa^2)$$
$$D_{6\tau}=3m^2+\kappa^2.$$

Apr\` es simplification, on obtient
$$m(Q_6)=\frac {\sqrt{3}}{48\pi^3}\sum'_{m,\kappa} (\frac {288\times 4 (m^2-3\kappa^2)}{(m^2+3\kappa^2)^3} +\frac {288\times 48\kappa^2}{(m^2+12\kappa^2)^3}-\frac {288\times 4^2\kappa^2}{(3m^2+4\kappa^2)^3}$$
$$-\frac {288}{(m^2+12\kappa^2)^2}+\frac {288}{(3m^2+4\kappa^2)^2})$$

De m\^ {e}me, pour $k=-36$, on a $\tau=\sqrt{-3}/3$ et 

$$D_{\tau}=\frac {1}{3} (m^2+3\kappa^2)$$
$$D_{2\tau}=\frac {1}{3} (4m^2+3\kappa^2)$$

$$D_{3\tau}=3m^2+\kappa^2)$$
$$D_{6\tau}=12m^2+\kappa^2.$$

Apr\` es simplification, on obtient
$$m(Q_{-36})=\frac {\sqrt{3}}{24\pi^3}\sum'_{m,\kappa} (\frac {72 (3\kappa^2-m^2)}{(m^2+3\kappa^2)^3} -\frac {288\times 12\kappa^2}{(4m^2+3\kappa^2)^3}+\frac {288\times 4\kappa^2}{(12m^2+\kappa^2)^3}$$
$$+\frac {288}{(4m^2+3\kappa^2)^2}-\frac {288}{(12m^2+\kappa^2)^2})$$
Comme
$$m(Q_0)=\frac {\sqrt {3}}{12\times 8\pi^3}4^3\times 3^2\sum'_{m,\kappa}\frac {\kappa^2-3m^2}{(\kappa^2+3m^2)^3}$$
on voit ais\'ement que la relation
$$2m(Q_{-36}-4m(Q_6)-m(Q_0)=0$$
n'est autre que la relation
$$\sum'_{m,\kappa}\frac {m^2-3\kappa^2}{(m^2+3\kappa^2)^3}=\sum'_{m,\kappa}(\frac {4m^2-3\kappa^2}{(4m^2+3\kappa^2)^3}-\frac {12m^2-\kappa^2}{(12m^2+\kappa^2)^3},$$
qui n'est autre que
 
$$L_{R'}(\phi_{R'},3)=L_R(\phi_R,3).$$

\qed

\textbf {Remerciements} \quad Je remercie vivement D. Boyd pour ses encouragements, ses conseils, ses v\'erifications, en un mot pour tout l'int\'er\^ {e}t qu'il a port\'e \` a ce travail.

\end{document}